\documentclass{amsart}

\newtheorem{theorem}{Theorem}[section]

\newtheorem{prop}[theorem]{Proposition}
\theoremstyle{definition}
\newtheorem{defn}[theorem]{Definition}
\newtheorem{examp}[theorem]{Example}
\newtheorem{exerc}[theorem]{Exercise}
\newtheorem{rem}[theorem]{Remark}
\numberwithin{equation}{section}

\newcommand{\To}{\longrightarrow}

\newcommand{\Gg}{\mathcal{G}}
\newcommand{\Jj}{\mathcal{J}}

\newcommand{\R}{\mathbb{R}}

\newcommand{\C}{\mathbb{C}}
\newcommand{\CP}{\mathbb{CP}}

\newcommand{\T}{{\mathbb T}}
\newcommand{\Z}{{\mathbb Z}}

\newcommand{\QED}{\hfill {\bf Q.E.D.} \medskip}

\newcommand{\al}{{\alpha}}

\newcommand{\Om}{{\Omega}}
\newcommand{\om}{{\omega}}
\newcommand{\de}{{\delta}}
\newcommand{\De}{{\Delta}}
\newcommand{\ga}{{\gamma}}

\newcommand{\la}{{\lambda}}

\newcommand{\Diff}{{\rm Diff}}
\newcommand{\Hess}{{\rm Hess}}

\newcommand{\disps}{\displaystyle}

\newcommand{\ov}{\overline}
\newcommand{\op}{{\overline{\partial}}}
\newcommand{\p}{{\partial}}

\newcommand{\tx}{\tilde{x}}
\newcommand{\tvphi}{\tilde{\varphi}}

\newcommand{\cp}{{{\CP}\,\!^2}}
\newcommand{\ocp}{{\ov{\CP}\,\!^2}}

\begin{document}

\title[K\"ahler geometry of toric manifolds]
{K\"ahler geometry of toric manifolds in symplectic coordinates}

\author{Miguel Abreu}

\address{Departamento de Matem\'{a}tica, Instituto Superior T\'{e}cnico,
Av.Rovisco Pais, 1049-001 Lisboa, Portugal}
\email{mabreu@math.ist.utl.pt}
\thanks{Partially supported by FCT grant PCEX/C/MAT/44/96 and by
PRAXIS XXI through the Research Units Pluriannual Funding Program.}



\date{\today}



\begin{abstract}
A theorem of Delzant states that any symplectic manifold 
$(M,\om)$ of dimension $2n$, equipped with an effective Hamiltonian action 
of the standard $n$-torus $\T^n = \R^{n}/2\pi\Z^n$, 
is a smooth projective toric variety completely 
determined (as a Hamiltonian $\T^n$-space) by the image of the moment 
map $\phi:M\to\R^n$, a convex polytope $P=\phi(M)\subset\R^n$. In 
this paper we show, using symplectic (action-angle) coordinates 
on $P\times \T^n$, how all $\om$-compatible toric complex structures on 
$M$ can be effectively parametrized by smooth functions on $P$. We also 
discuss some topics suited for application of this symplectic coordinates 
approach to K\"ahler toric geometry, namely: explicit construction of 
extremal K\"ahler metrics, spectral properties of toric manifolds and 
combinatorics of polytopes.
\end{abstract}

\maketitle

\section{Introduction} \label{sec:intro}

K\"ahler geometry can be thought of as a ``compatible'' intersection 
of complex and symplectic geometries. Indeed, the triple 
$(M^{2n},J,\om)$, with $2n$ the {\bf real} dimension of $M$, is a 
{\bf K\"ahler manifold} if
\begin{itemize}
\item[(i)] $(M^{2n}, J)$ is a {\bf complex manifold}, i.e. the 
automorphism $J:TM\to TM$, $J^2 = -I$, is an integrable complex 
structure;
\item[(ii)] $(M^{2n},\om)$ is a {\bf symplectic manifold}, i.e. the 
$2$-form $\om$ is closed and non-degenerate;
\item[(iii)] $J$ and $\om$ are {\bf compatible} in the sense that the 
bilinear form $\om (\cdot, J\cdot )$ is a {\bf Riemannian metric}, 
i.e. symmetric and positive definite.
\end{itemize}

There are two natural types of local coordinates on a K\"ahler 
manifold: complex (holomorphic) and symplectic (Darboux) coordinates. 
In almost all differential geometric questions on K\"ahler manifolds,
like local Riemannian K\"ahler geometry or existence of 
K\"ahler-Einstein metrics, the complex point of view plays the 
dominant role. There are several reasons for that, one being that in 
complex coordinates all compatible symplectic forms can be effectively 
parametrized  using functions. In fact, if one fixes a complex 
structure $J_{0}$ and local holomorphic coordinates $(z_{1},\ldots, 
z_{n})$, any compatible symplectic form $\om$ is locally of the form
\begin{equation} \label{localsymp}
\om = 2i\p\op f,\ f\equiv\mbox{local smooth real function.}
\end{equation}
Moreover, if one is given a fixed compatible symplectic form 
$\om_{0}$ on $M$, defining a cohomology class $\Om \in H^2 (M,\R)$, 
any other compatible symplectic form $\om_{1}$ in the same class $\Om$ 
is globally of the form
\begin{equation} \label{globalsymp}
\om_{1} = \om_{0} + 2i\p\op f, \ f\in C^\infty (M).
\end{equation}

Because the family $\om_{t} = \om_{0} + t (\om_{1} - \om_{0}),\ t\in 
[0,1]$, is an isotopy of symplectic forms in the same cohomology 
class $\Om$, Moser's theorem~\cite{mose} gives a family of 
diffeomorphisms $\varphi_{t}:M\to M,\ t\in [0,1]$, such that 
$\varphi_{t}^\ast (\om_{t}) = \om_{0}$. Hence $\varphi_{1}$ is a
K\"ahler isomorphism between $(M,J_{0},\om_{1})$ and $(M,J_{1},\om_{0})$,
where $J_{1} = (\varphi_{1})_{\ast}^{-1} \circ J_{0} \circ 
(\varphi_{1})_{\ast}$. This means that {\it fixing the complex structure 
$J_{0}$ and varying the compatible symplectic form $\om$ in a fixed 
cohomology class is equivalent to fixing the symplectic form $\om_{0}$ 
and varying the compatible complex structure $J$ in a fixed 
diffeomorphism class}. This simple fact makes it reasonable, at least 
in principle, to try to do K\"ahler geometry starting from a 
symplectic coordinate chart, i.e. Darboux coordinates 
$(x_{1},\ldots,x_{n},y_{1},\ldots,y_{n})$ on which $\om_{0}$ has the 
standard form
$$ \om_{0} = \sum_{j} dx_{j}\wedge dy_{j}\ .$$
Unfortunately this does not take us very far in general, since one 
does not know of any effective way of parametrizing (locally 
or globally) compatible complex structures.

Nevertheless, the main purpose of this paper is to illustrate how
\begin{center}
{\bf symplectic coordinates can be useful in K\"ahler geometry,}
\end{center}
more precisely in K\"ahler toric geometry. As is explained in 
\S\ref{ssec:def}, a K\"ahler toric manifold $M$ can be described 
either as a compactification of a complex torus $\T_{\C}^n = \C^n / 
2\pi i\Z^n$ (complex point of view), or as a compactification of 
$P^{\circ}\times \T^n$ (symplectic point of view), where 
$P^{\circ}\subset\R^n$ is the interior of a certain convex polytope 
(moment polytope), $\T^n$ is a real torus $\R^n / 2\pi\Z^n$, and the 
symplectic form is the standard $\om_{0} = \sum_{i} dx_{i}\wedge dy_{i}$ 
(the $x$'s and $y$'s being linear coordinates on $P^\circ\subset\R^n$ 
and $\T^n$, respectively). In \S\ref{ssec:set-up}, which in part is based 
on Guillemin's paper~\cite{gui2}, we will see how to change from 
one point of view to the other and back again by an appropriate 
Legendre transform and its dual, and also how all $\om_{0}$-compatible toric 
complex structures can be effectively parametrized using smooth 
functions on $P$ (more precisely, through their ``Hessian''). Then, in 
\S\ref{ssec:can}, we present Guillemin's explicit formula for a 
``canonical'' $\om_{0}$-compatible toric complex structure $J_{0}$ on 
$M$, expressed only in terms of combinatorial data on $P$. The 
analogue of~(\ref{globalsymp}) for $\om_{0}$-compatible toric complex 
structures is presented in \S\ref{ssec:general} and proved in the Appendix.
Finally, in section~\ref{sec:applications}, we discuss some applications and 
open problems related to this symplectic approach to K\"ahler toric
geometry. Specifically, we talk about:
\begin{itemize}
\item[-] explicit constructions of extremal K\"ahler metrics in 
\S\ref{ssec:extremal};
\item[-] spectral properties of toric manifolds in 
\S\ref{ssec:spectral};
\item[-] relation between K\"ahler toric geometry and combinatorics of 
polytopes in \S\ref{ssec:combinatorics}.
\end{itemize}
These topics, chosen according to the author's current knowledge and previous 
work, are obviously not exhaustive. We believe that other 
differential-geometric questions can be successfully approached through 
this ``symplectic coordinates'' point of view.

My interest in this subject originated from Guillemin's 
paper~\cite{gui2}. The reading of recent papers by 
Donaldson~(\cite{don1} and~\cite{don2}) and Hitchin~(\cite{hit1} 
and~\cite{hit2}), where this ``duality'' between changes in complex 
structure and changes in symplectic structure is also present (in the 
context, respectively, of symplectic/complex quotients of infinite 
dimensional spaces and the Strominger-Yau-Zaslow approach to mirror 
symmetry), contributed to the clarification in my mind of some of the 
ideas presented in this paper.

\section{K\"ahler toric manifolds} \label{sec:toric}

In this section we present a description of toric manifolds and their 
invariant K\"ahler metrics. The point of view is that of a symplectic 
geometer trying to translate into its own natural context and 
language some well-known facts from K\"ahler geometry. Good 
references for most of what will be presented in this section are 
Guillemin's paper~\cite{gui2}, book~\cite{gui1} and some of the 
bibliography listed there. 

\subsection{Definition and ``canonical'' examples} \label{ssec:def}

\begin{defn} \label{def:toric}
A {\bf K\"ahler toric manifold} is a closed connected 
$2n$-dimensional K\"ahler manifold $(M^{2n},\om,J)$ equipped with an 
effective Hamiltonian holomorphic action
$$\tau:\T^n \rightarrow \Diff(M^{2n},\om,J)$$
of the standard (real) $n$-torus $\T^n$.
\end{defn}

A few comments are in order regarding this definition. The category 
where it fits is, of course, K\"ahler geometry and not just complex 
geometry or symplectic geometry. As we recall below, almost all closed 
smooth complex toric varieties and all closed symplectic toric 
manifolds are K\"ahler. However, in the complex or symplectic 
categories, the particular form of a K\"ahler metric determined by a 
compatible symplectic form or complex structure is secondary and not 
part of the original data, although the fact that one exists is 
important and often used. In this paper, in particular 
Definition~\ref{def:toric}, the object of study is the K\"ahler metric 
involving both the complex and symplectic structures. Two K\"ahler 
toric manifolds are isomorphic if they are equivariantly K\"ahler 
isomorphic, which implies in particular that they have to be 
isometric as Riemannian manifolds. For example, while in the purely 
complex or symplectic (with fixed volume) worlds there is only one 
toric manifold of dimension $2$, the standard sphere in $\R^3$, in the 
K\"ahler world different $S^1$-invariant Riemannian metrics on $S^2$ 
determine different K\"ahler toric manifolds of dimension $2$ (see 
\S\ref{ssec:spectral}). Related to this, see also Remark~\ref{rmk:Pdet}.

We now recall the usual complex and symplectic definitions of a 
toric manifold, not only for comparing with 
Definition~\ref{def:toric}, but also to introduce some of the purely 
complex and symplectic ingredients we will need.

From the complex geometry point of view, see~\cite{dani} or~\cite{fult}, 
a toric variety is a normal variety $X$ that contains a complex torus 
$\T_{\C}^n = \C^n / 2\pi i \Z^n$ as a dense open subset, together with 
an holomorphic action 
$$\T_{\C}^n\times X \to X$$
of $\T_{\C}^n$ on $X$ that extends the natural action of $\T_{\C}^n$ 
on itself. Hence, a toric variety is a suitable compactification of 
the complex torus $\T_{\C}^n$. Each particular compactification is 
determined by a combinatorial gadget called a ``fan''. Here we will 
not discuss the details of this complex point of view, although we will 
make some use of the holomorphic coordinate description of the open 
dense $\T_{\C}^n$ orbit as $\C^n / 2\pi i \Z^n$ (which, as will be seen 
in the Appendix, is always present on a K\"ahler toric manifold).

Not all closed smooth complex toric varieties are projective 
(example on p.71 of~\cite{fult}), but they are all dominated 
birationally by a smooth projective toric variety (exercise on p.72 
of~\cite{fult}). Any smooth projective toric variety is K\"ahler and 
an example of a K\"ahler toric manifold (the action of the real torus
$$\T^n = i\R^n / 2\pi i \Z^n \subset \C^n / 2\pi i \Z^n = \T_{\C}^n$$
is Hamiltonian with respect to the K\"ahler form). Moreover the 
converse is also true, i.e. any K\"ahler toric manifold is 
equivariantly biholomorphic (but not necessarily K\"ahler isomorphic)
to a smooth projective toric variety (see Remark~\ref{rmk:Pdet} and 
Proposition~\ref{prop:uniqueJ} in the Appendix).

From the symplectic geometry point of view, a symplectic toric 
manifold is a symplectic manifold $(M^{2n},\om)$ equipped with an 
effective Hamiltonian action
$$\tau:\T^n \rightarrow \Diff(M^{2n},\om)$$
of the standard (real) $n$-torus $\T^n = \R^n / 2\pi\Z^n$. Associated 
to $(M,\om,\tau)$ there is a {\bf moment map} $\phi:M\to\R^n$, whose 
image $P=\phi(M)\subset\R^n$ is a convex polytope. Denoting by 
$P^\circ$ the interior of $P$, we have that the set $M^\circ = 
\phi^{-1}(P^\circ)$ is an open dense subset of $M$ consisting of the 
points where the $\T^n$-action is free. $M^\circ$ is symplectomorphic 
to $P^\circ \times \T^n \subset \R^n \times \T^n$ (with its standard 
linear symplectic structure induced from $\R^{2n}$). Hence, a 
symplectic toric manifold can be thought of as a suitable 
compactification of an $n$-parameter family of Lagrangian tori $\T^n$, 
the compactification being determined by the combinatorics of a convex 
polytope $P\subset \R^n$.

Not every convex polytope in $\R^n$ is the moment polytope of some triple
$(M,\om,\tau)$. The following definition characterizes the ones that are.
\begin{defn} \label{def:delzant} 
A convex polytope $P$ in $\R^n$ is {\bf Delzant} if:
\begin{description}
\item[(1)] there are $n$ edges meeting at each vertex $p$;
\item[(2)] the edges meeting at the vertex $p$ are rational, i.e. each edge
is of the form $p + tv_i,\ 0\leq t\leq \infty,\ {\rm where}\ v_i\in\Z^n$;
\item[(3)] the $v_1, \ldots, v_n$ in (2) can be chosen to be a basis of
$\Z^n$.
\end{description}
\end{defn}
In \cite{delz} Delzant associates to every Delzant polytope $P\subset\R^n$
a closed connected symplectic manifold $(M_P,\om_{P})$ of dimension 
$2n$, together with a Hamiltonian $\T^n$-action 
$$\tau_{P}:\T^n \rightarrow \Diff(M_{P},\om_{P})$$
with moment map 
$\phi_{P}:M_{P}\to\R^n$, such that the image $\phi_{P}(M)$ is 
precisely $P$. Moreover, he shows that this is a bijective
correspondence. More precisely, he proves the following:
\begin{theorem} \label{th:delzant}
Let $(M,\om,\tau)$ be a compact, connected, $2n$-dimensional
Hamiltonian $\T^n$-space, on which the action of $\T^n$ is effective with
moment map $\phi : M\rightarrow\R^n$. Then the image $P$ of $\phi$ is a
Delzant polytope, and $(M,\om,\tau)$ is isomorphic as a Hamiltonian 
$\T^n$-space to $(M_{P},\om_{P},\tau_{P})$.
\end{theorem}

It follows from Delzant's construction that $(M_P,\om_{P},\tau_{P})$ 
(or, because of the
above theorem, any effective Hamiltonian $\T^n$-space with the same
moment polytope) is equipped with a ``canonical'' $\T^n$-invariant complex
structure $J_{P}$ compatible with the symplectic form $\om_{P}$.
In other words, for each Delzant polytope $P\subset\R^n$, the quadruple 
$(M_{P},\om_{P},J_{P},\tau_{P})$ is an example of a K\"ahler toric manifold.
These are the {\bf ``canonical'' examples} of the title of this subsection.
Delzant's construction also shows that the $\T^n$-action can 
be complexified to a holomorphic $\T_{\C}^n$-action on $(M_{P},J_{P})$,
giving $M_{P}$ the structure of a smooth projective toric variety.
All this ``canonical'' examples are actually K\"ahler projective, in 
the sense that not only $J_{P}$ but also $\om_{P}$ is induced from the 
standard Fubini-Study one on complex projective space.
\begin{rem} \label{rmk:Pdet} For any K\"ahler toric manifold 
$(M,\om,J,\tau)$, with moment polytope $P\subset\R^n$, we have that:
\begin{itemize}
\item[(i)] $(M,\om,\tau)$ is equivariantly symplectomorphic to the 
``canonical'' $(M_{P},\om_{P},\tau_{P})$ (by Theorem~\ref{th:delzant});
\item[(ii)] $(M,J,\tau)$ is equivariantly biholomorphic to the 
``canonical'' $(M_{P},J_{P},\tau_{P})$ (see 
Proposition~\ref{prop:uniqueJ} in the Appendix).
\end{itemize}
Hence, the polytope $P$ completely determines the symplectic and 
complex structures of $(M,\om,J,\tau)$, but only if we consider them 
separately. The toric K\"ahler metric $\om(\cdot, J\,\cdot)$,
obtained by combining the two, is {\bf not} determined by $P$.
\end{rem}

\subsection{K\"ahler metrics in symplectic coordinates: set-up} 
\label{ssec:set-up}

Here we explain the set-up that will be used to describe 
invariant K\"ahler metrics on a toric manifold in symplectic 
coordinates. We first recall the usual set-up in complex 
coordinates, then present the one in symplectic coordinates and 
finally justify it by showing how the moment map (or Legendre 
transform) is the explicit translator from complex to symplectic.

Let $(M^{2n},\om,J,\tau)$, $\tau:\T^n \to \Diff (M,\om,J)$, be a 
K\"ahler toric manifold. We again denote by $M^\circ$ the open dense 
subset of $M$ defined by
$$M^\circ = \left\{ p\in M : \T^n\mbox{-action is free at $p$}\right\}.$$
We can describe $M^\circ$ in complex (holomorphic) coordinates as
$$M^\circ \cong \C^n / 2\pi i\Z^n = \R^n \times i \T^n =
\left\{ u+iv : u\in\R^n,\ v\in \R^n / \Z^n\right\}$$
(see the proof of Proposition~\ref{prop:uniqueJ} in the Appendix).
In this $z=u+iv$ coordinates, the $\T^n$-action is given by
$$t\cdot (u+iv) = u + i(v+t),\ t\in\T^n,$$
and the complex structure $J$, which is just multiplication by $i$, is 
given by the standard
$$ J = \left[ \begin{array}{ccc}
0 & \vdots & -I \\
\cdots & \cdots & \cdots \\
I & \vdots & 0
\end{array}\right]$$
where $I$ denotes the $(n\times n)$ identity matrix. The symplectic 
K\"ahler form $\om$ is given by a potential $f\in C^\infty 
(M^\circ)$ through the usual relation $\om = 2i\p\op f$. Since $\om$ 
is invariant by the $\T^n$-action, the potential $f$ depends only on 
the $u$ coordinates: $f = f(u) \in C^\infty(\R^n)$. Hence, the matrix 
that represents the skew-symmetric bilinear form $\om$ in this $(u,v)$ 
coordinates is of the form
$$\left[ \begin{array}{ccc}
0 & \vdots & F \\
\cdots & \cdots & \cdots \\
-F & \vdots & 0
\end{array}\right]$$
with $F = \Hess_{u}(f) \equiv$ Hessian of $f$ in the $u$ coordinates:
$$F = \left[ f_{jk} \right]_{j,k=1}^{n,n},\ 
f_{jk} = \frac{\p^2 f}{\p u_{j} \p u_{k}},\ 1\leq j,k\leq n\ .$$
Note that the Riemannian K\"ahler metric $\om(\cdot, J\cdot)$ is then 
given in matrix form by
\begin{equation} \label{metricomp}
\left[ \begin{array}{ccc}
F & \vdots & 0\\
\cdots & \cdots & \cdots \\
0 & \vdots & F
\end{array}\right]
\end{equation}

Not every function $f\in C^\infty(\R^n)$ is the K\"ahler potential of 
a symplectic form $\om$. Because~(\ref{metricomp}) defines a 
Riemannian metric, we have that $F = \Hess_{u}(f)$ is positive definite.
This means that $f$ is strictly convex. Moreover, the fact that 
$\om$ and the metric compactify smoothly to give a K\"ahler form 
and metric on $M$, puts restrictions on the behaviour of $f$ at 
infinity. However we will not have to worry with what those 
restrictions are since, as we will see in \S~\ref{ssec:can} 
and \S~\ref{ssec:general}, they become simple and explicit in the 
symplectic set-up that we are interested in and will work 
with.

In symplectic (or action-angle) coordinates, $M^\circ$ can be 
described as
$$M^\circ \cong P^\circ \times \T^n = \left\{ (x,y) :
x\in P^\circ\subset\R^n, y\in\R^n / \Z^n \right\}$$
where $P^\circ$ is the interior of a Delzant polytope $P\subset\R^n$.
In these $(x,y)$ coordinates, the $\T^n$-action is given by
$$t\cdot (x,y) = (x,y+t),\ t\in\T^n,$$
and the symplectic form $\om$ is the standard $\om = dx\wedge dy
= \sum_{j}dx_{j}\wedge dy_{j}$, which in matrix form is
$$\left[ \begin{array}{ccc}
0 & \vdots & I \\
\cdots & \cdots & \cdots \\
-I & \vdots & 0
\end{array}\right]$$
The interesting part now is how one describes the complex K\"ahler 
structure $J$. It turns out, and we will see below why, that $J$ is 
given by a ``potential'' $g=g(x)\in C^\infty (P^\circ)$ through the 
relation $J =\ \mbox{``Hessian'' of}\ g$. More precisely, the matrix 
that represents the complex structure $J$ in these $(x,y)$ coordinates 
is of the form
\begin{equation} \label{Jsymp}
J = \left[ \begin{array}{ccc}
0 & \vdots & -G^{-1} \\
\cdots & \cdots & \cdots \\
G & \vdots & 0
\end{array}\right]
\end{equation}
with $G = \Hess_{x}(g) \equiv$ Hessian of $g$ in the $x$ coordinates:
$$G = \left[ g_{jk} \right]_{j,k=1}^{n,n},\ 
g_{jk} = \frac{\p^2 g}{\p x_{j} \p x_{k}},\ 1\leq j,k\leq n\ .$$
Note that the Riemannian K\"ahler metric $\om(\cdot, J\cdot)$ is now 
given in matrix form by
\begin{equation} \label{metricsymp}
\left[ \begin{array}{ccc}
G & \vdots & 0\\
\cdots & \cdots & \cdots \\
0 & \vdots & G^{-1}
\end{array}\right]
\end{equation}

Let us address the question of how to change from one set of 
coordinates to the other and vice-versa. In complex coordinates, 
the moment map $\phi$ of the $\T^n$-action with respect to $\om$ is 
given by the Legendre transform
$$\phi(u,v) = \frac{\p f}{\p u}\ .$$
The restriction $\phi(u,0),\ u\in\R^n,$ is a diffeomorphism of $\R^n$ 
onto $P^\circ$. Hence the map
\begin{equation} \label{hol-symp}
x = \frac{\p f}{\p u},\ y = v,
\end{equation}
gives us a diffeomorphism from $\R^n\times\T^n$ to $P^\circ\times\T^n$.
Moreover, the symplectic form transforms under this change of 
coordinates from $\om = 2i\p\op f$ to the standard $\om = dx \wedge 
dy$.

The complex structure $J$, standard in the $(u,v)$ coordinates, 
clearly becomes of the form~(\ref{Jsymp}) where $G$ at the point 
$x=\frac{\p f}{\p u}$ is equal to the inverse of $F$ at the point $u$.
We can now argue in two different ways to show that $G$ is of the 
form $\Hess_{x}(g)$ for some $g\in C^\infty (P^\circ)$. One is to 
quote some general facts about the Legendre transform (see pp.118-126 
of~\cite{gui1}). The other is to check that the integrability 
condition for $J$ translates into the following compatibility 
conditions for $G = \left[ g_{jk} \right]$:
$$g_{jk,l} \equiv \frac{\p g_{jk}}{\p x_{l}} = \frac{\p g_{jl}}{\p x_{k}}
\equiv g_{jl,k},\  1\leq j,k,l \leq n.$$
Hence we do have a function $g\in C^\infty (P^\circ)$ such that 
$G = \Hess_{x}(g)$. Since $G$ at the point $x=\frac{\p f}{\p u}$ is equal 
to the inverse of $F$ at the point $u$, we also have that the map
\begin{equation} \label{symp-hol}
u = \frac{\p g}{\p x},\ v = y,
\end{equation}
is (up to a constant) the inverse of~(\ref{hol-symp}). In other 
words, $f$ and $g$ are (up to a linear factor) Legendre dual to each 
other:
\begin{equation}\label{duality}
f(u)+g(x) = \sum_{j}\frac{\p f}{\p u_{j}}\cdot\frac{\p g}{\p x_{j}},\ 
\mbox{at}\ x=\frac{\p f}{\p u}\ \mbox{or}\ u=\frac{\p g}{\p x}.
\end{equation}

\subsection{K\"ahler metrics in symplectic coordinates: ``canonical'' 
examples} \label{ssec:can}

As was already mentioned in \S\ref{ssec:def}, a construction of 
Delzant~\cite{delz} associates to every Delzant polytope 
$P\subset\R^n$ a ``canonical'' K\"ahler toric manifold 
$(M_{P},\om_{P},J_{P},\tau_{P})$. In~\cite{gui2} Guillemin gives an explicit 
formula for the K\"ahler metric $\om_{P}(\cdot,J_{P}\,\cdot)$
in terms of combinatorial data on 
$P$. The purpose of this subsection is to present that formula in the 
set-up just explained.

A Delzant polytope $P$ can be described by a set of inequalities of 
the form
$\langle x,\mu_{r}\rangle \geq \la_{r},\ r=1,\ldots,d$, each 
$\mu_{r}$ being a primitive element of the lattice $\Z^n\subset\R^n$ 
and inward-pointing normal to the $r$-th $(n-1)$-dimensional face of 
$P$. Consider the affine functions $\ell_{r}:\R^n\to\R,\ 
r=1,\ldots,d$, defined by
\begin{equation} \label{ells} 
\ell_{r}(x) = \langle x,\mu_{r}\rangle - \la_{r}\ .
\end{equation}
Then $x\in P^\circ$ if and only if $\ell_{r}(x) > 0$ for all $r$, and 
hence the function $g_{P}:P^\circ\to\R$ defined by
\begin{equation} \label{gstandard}
g_{P}(x) = \frac{1}{2}\sum_{r=1}^{d} \ell_{r}(x) \log \ell_{r}(x)
\end{equation}
is smooth on $P^\circ$.
\begin{theorem}[Guillemin,\cite{gui2}] \label{th:Jcan} The ``canonical'' 
compatible toric complex structure $J_{P}$ on $(M_{P},\om_{P})$ is 
given in the $(x,y)$ symplectic coordinates of $M_{P}^\circ \cong 
P^\circ\times\T^n$ by
$$
J_{P} = \left[ \begin{array}{ccc}
0 & \vdots & - G_{P}^{-1} \\
\cdots & \cdots & \cdots \\
G_{P} & \vdots & 0
\end{array}\right]
$$
with $G_{P} = \Hess_{x}(g_{P}),\ g_{P}\equiv$ smooth function on
$P^{\circ}$ defined by~(\ref{gstandard}).
\end{theorem}
Note that, as it should be, $g_{P}$ is singular on the boundary of 
$P^\circ$ (and so are $G_{P}$ and $J_{P}$). A consequence of 
Guillemin's theorem is that the type of singular behaviour that 
$g_{P}$ has on $\p P^\circ$ is the right one to allow for the smooth 
compactification of $J_{P}$ (and also the associated Riemannian 
K\"ahler metric $\om_{P}(\cdot,J_{P}\,\cdot)$) to the whole $M_{P}$.

\begin{examp}[$2$-sphere] \label{ex:sphere}
Let $S^2\subset\R^3$ be the standard sphere of radius $1$, with 
symplectic form given by the standard area form with total area 
$4\pi$, and $S^1$-action given by rotation around an axis. The moment 
map $\phi:S^2\to\R$ is just projection to this axis of rotation, and 
so the moment polytope is given by $P=[-1,1]$. This is determined by 
the two affine functions
$$\ell_{-1}(x) = 1+x\ \mbox{and}\ \ell_{1}(x) = 1-x\ ,$$
which means that the function $g_{P}:P^\circ = (-1,1)\to\R$ has the 
form
$$g_{P}(x) = \frac{1}{2} \left[ (1+x)\log(1+x) + (1-x)\log(1-x) 
\right]\ .$$
Computing two derivatives on gets that
$$G_{P} = \left[ \frac{1}{1-x^2} \right]\ ,$$
and so, on $P^\circ\times S^1 = (-1,1)\times S^1$ with $(x,y)$ 
coordinates, we have that
$$\om_{P} = dx\wedge dy\,;\ 
J_{P} = \left[ \begin{array}{cc} 0 & - (1 - x^2) \\ 
\frac{\displaystyle 1}{\displaystyle 1 - x^2} &  0
\end{array}\right]\,;\ 
\om_{P}(\cdot,J_{P}\,\cdot) = \left[ \begin{array}{cc}
\frac{\displaystyle 1}{\displaystyle 1 - x^2} & 0 \\ 
0 & (1 - x^2) \end{array}\right].$$
As we will see in \S\ref{ssec:extremal}, and one can easily check, 
these correspond to the standard area form, complex structure and 
round metric on $S^2$.
\end{examp}

\begin{examp}[Projective plane] \label{ex:cp2}
Let $P\subset\R^2$ be the Delzant polytope defined by
$$\ell_{1}(x_{1},x_{2}) = 1+x_{1}\,,\ \ell_{2}(x_{1},x_{2}) = 1+x_{2}\ 
\mbox{and}\ \ell_{3}(x_{1},x_{2}) = 1 - x_{1} - x_{2}\ ,$$
an isosceles right triangle that is well-known to correspond to 
$\cp$ (see e.g. Chapter $1$ of~\cite{gui1}). We then have that
\begin{eqnarray*}
g_{P}(x_{1},x_{2}) & = & \frac{1}{2} \left[ (1+x_{1})\log(1+x_{1}) + 
                         (1+x_{2})\log(1+x_{2}) + \right. \\
				   &   & \ \ \ \,\left. (1-x_{1}-x_{2})\log(1-x_{1}-x_{2})
				                 \right]
\end{eqnarray*}
and
$$ G_{P}(x_{1},x_{2}) = \Hess_{x}(g_{P}) = \frac{1}{2(1-x_{1}-x_{2})}
\left[ \begin{array}{cc} \frac{\displaystyle 2-x_{2}}{\displaystyle 1+x_{1}} 
& 1 \\ 1 &  \frac{\displaystyle 2-x_{1}}{\displaystyle 1+x_{2}}
\end{array}\right]\ .$$
One can check, either by tracing down Delzant's construction or using
the curvature computations of \S\ref{ssec:extremal}, that
this corresponds to the standard Fubini-Study K\"ahler 
metric on $\cp$ with area $6\pi$ on ${{\CP}\,\!^1}$.
\end{examp}

\subsection{K\"ahler metrics in symplectic coordinates: general 
description} \label{ssec:general}

Our goal in this subsection is to describe all possible compatible 
toric complex structures $J$ on the symplectic toric manifold 
$(M_{P},\om_{P},\tau_{P})$ associated to a Delzant polytope $P\subset\R^n$. If 
possible, we would like this description to be similar to the one 
for compatible symplectic forms given by~(\ref{globalsymp}) .

In light of the set-up of \S\ref{ssec:set-up}, in 
particular~(\ref{Jsymp}), any such $J$ is determined in the $(x,y)$ 
symplectic coordinates of $M^\circ \cong P^\circ\times\T^n$ by a 
``potential'' $g=g(x)\in C^\infty(P^\circ)$. Because of the 
form~(\ref{metricsymp}) of the Riemannian K\"ahler metric 
$\om_{P}(\cdot,J\cdot)$, we have that $G=\Hess_{x}(g)$ is positive 
definite and so $g$ has to be strictly convex on $P^\circ$.

To understand what the behaviour of $g$ has to be near the boundary 
of $P^\circ$, so that $\om_{P}$ and $J$ compactify smoothly on 
$M_{P}$, it is useful to look more carefully at the standard $g_{P}$, 
given by~(\ref{gstandard}), and the corresponding $G_{P}$. Explicit 
calculations show easily that, although $G_{P}$ is singular on the 
boundary of $P^\circ$, $G_{P}^{-1}$ is smooth on the whole $P$ and its 
determinant has the form
$$\det(G_{P}^{-1}) = \de_{P}(x)\cdot\prod_{r=1}^{d} \ell_{r}(x)\ ,$$
where the $\ell_{r}$'s are defined by~(\ref{ells}) and $\de_{P}$ is a 
smooth and strictly positive function on the whole $P$.
The geometric interpretation of this formula is clear. As one reaches 
the $r$-th $(n-1)$-dimensional face of $P$, the positive definite 
matrix $G_{P}^{-1}$ acquires a kernel that is generated by the 
normal $\mu_{r}$. Because of the nondegeneracy properties 
of a Delzant polytope (see Definition~\ref{def:delzant}), this means for 
example that at any vertex $p$ of $P$ we have $G_{P}^{-1}\equiv$ zero 
matrix.

This type of singular behaviour of $g_{P}$, and degenerate behaviour 
of $G_{P}^{-1}$, has to be present in any $g$ and $G^{-1}$ that 
correspond to honest smooth symplectic and complex structures on the 
whole $M_{P}$. The following theorem, which describes all possible 
compatible toric $J$'s in a way very similar to~(\ref{globalsymp}), 
is just a consequence of that. A detailed proof is given in the 
Appendix.
\begin{theorem} \label{th:Jgen}
Let $(M_{P},\om_{P},\tau_{P})$ be the toric symplectic manifold associated to 
a Delzant polytope $P\subset\R^n$, and $J$ any compatible toric 
complex structure. Then $J$ is determined, using~(\ref{Jsymp}),
by a ``potential'' $g\in C^\infty(P^\circ)$ of the form
\begin{equation} \label{ggeneral}
g = g_{P} + h\ ,
\end{equation}
where $g_{P}$ is given by~(\ref{gstandard}), $h$ is smooth on the 
whole $P$, and the matrix $G = \Hess_{x}(g)$ is positive definite on
$P^\circ$ and has determinant of the form
$$\det(G) = \left[ \de(x)\cdot\prod_{r=1}^{d} \ell_{r}(x) 
\right]^{-1}\ ,$$
with $\de$ being a smooth and strictly positive function on the whole $P$.

Conversely, any such $g$ determines a compatible toric complex 
structure $J$ on $(M_{P},\om_{P})$, which in the $(x,y)$ symplectic 
coordinates of $M_{P}^\circ \cong P^\circ\times\T^n$ has the 
form~(\ref{Jsymp}).
\end{theorem}

\section{Applications and some interesting problems} \label{sec:applications}

In this section we discuss some applications of the symplectic approach to
K\"ahler geometry of toric manifolds, presented in \S\ref{sec:toric},
and also state some open problems where we think this approach should be 
useful. These applications and problems, chosen solely according to the
author's current knowledge and previous work, deal with the following 
differential-geometric issues and tools: scalar 
curvature~(\S\ref{ssec:extremal}), spectrum of the 
Laplacian~(\S\ref{ssec:spectral}) and differential
forms~(\S\ref{ssec:combinatorics}).

\subsection{Explicit construction of extremal K\"ahler metrics} 
\label{ssec:extremal}

In~\cite{cal1} and~\cite{cal2}, Calabi introduced the notion of {\bf extremal}
K\"ahler metrics. These are defined, for a fixed closed complex manifold
$(M,J_0)$, as critical points of the square of the $L^2$-norm of the scalar
curvature, considered as a functional on the space of all symplectic K\"ahler
forms $\om$ in a fixed K\"ahler class $\Om\in H^2(M,\R)$. The extremal
Euler-Lagrange equation is equivalent to the gradient of the scalar curvature
being an holomorphic vector field (see~\cite{cal1}), and so these metrics
generalize constant scalar curvature K\"ahler metrics. Moreover, Calabi
showed in~\cite{cal2} that extremal K\"ahler metrics are always invariant under
a maximal compact subgroup of the group of holomorphic transformations of
$(M,J_0)$. Hence, on a complex toric manifold, extremal K\"ahler metrics are
automatically toric K\"ahler metrics, and one should be able to write them
down using the general description given by Theorem~\ref{th:Jgen}.

An attempt in this direction was made in~\cite{abr}, and we now briefly
summarize what is obtained there. A Riemannian K\"ahler metric, given in 
holomorphic coordinates $(u,v)$ by~(\ref{metricomp}), has scalar curvature
$S$ given by
\begin{equation} \label{scalarcomp}
S = - \frac{1}{2} \sum_{j,k} f^{jk}\, \frac{\p^2 \log (\det F)}
{\p u_j \p u_k}\, ,
\end{equation}
where the $f^{jk},\ 1\leq j,k\leq n$, denote the entries of the inverse of the
matrix $F = \Hess_u (f)$ (see~\S\ref{ssec:set-up}). Using~(\ref{hol-symp})
this can be written in symplectic $(x,y)$ coordinates as
\begin{equation} \label{scalarsymp1}
S = - \frac{1}{2} \sum_{j,k} \frac{\p}{\p x_j}
\left( g^{jk}\, \frac{\p \log (\det G)}{\p x_k} \right)\,,
\end{equation}
which after some algebraic manipulations becomes the more compact
\begin{equation} \label{scalarsymp2}
S = - \frac{1}{2} \sum_{j,k} \frac{\p^2 g^{jk}}{\p x_j \p x_k}\,, 
\end{equation}
where now the $g^{jk},\ 1\leq j,k\leq n$, are the entries of the inverse 
of the matrix $G = \Hess_x (g)$. The Euler-Lagrange equation defining an
extremal K\"ahler metric can be shown to be equivalent to
\begin{equation} \label{extremalsymp}
\frac{\p S}{\p x_j} \equiv\ \mbox{constant},\ j=1,\ldots,n,
\end{equation}
i.e. {\it the metric is extremal if and only if its scalar curvature $S$ is an
affine function of $x$}.

\begin{examp}
For the ``canonical" toric metric on $S^2$, presented in 
Example~\ref{ex:sphere}, we have that
$$G_P^{-1} = \left[ 1-x^2\right]$$
and so
$$S_P (x) = -\frac{1}{2} (1-x^2)''\equiv 1.$$

For the ``canonical" toric metric on $\cp$, presented in 
Example~\ref{ex:cp2}, we have that
$$ G_{P}^{-1}(x_{1},x_{2}) = \frac{2}{3}
\left[ \begin{array}{cc} (2-x_1)(1+x_1) & -(1+x_1)(1+x_2) \\
  & \\ -(1+x_1)(1+x_2) &  (2-x_2)(1+x_2)
\end{array}\right]$$
and so
$$S_P (x_1,x_2) \equiv -\frac{1}{3}\left[ -2-1-1-2\right] = 2\ .$$

In both these cases, the fact that the metrics have constant scalar curvature
makes them (trivial) examples of extremal K\"ahler metrics and can be used
to show that they correspond to standard metrics on $S^2$ and $\cp$.
\end{examp}

In~\cite{cal1}, Calabi constructed families of extremal K\"ahler metrics
of non-constant scalar curvature. In~\cite{abr}, Calabi's simplest family
(on $\cp\#\,\ocp$) is written down very simply and explicitly in 
symplectic coordinates, providing an example of the effectiveness of the 
parametrization given by Theorem~\ref{th:Jgen}. The following exercise 
describes a particular metric in that family.

\begin{exerc}[Blow-up of $\cp$] \label{ex:blow-up}
Let $P\subset\R^2$ be the Delzant polytope defined by
$$\ell_{1}(x_{1},x_{2}) = 1+x_{1}\,,\ \ell_{2}(x_{1},x_{2}) = 1+x_{2}\,,\  
\ell_{3}(x_{1},x_{2}) = 1 - x_{1} - x_{2}$$
and
$$\ell_{-3}(x_{1},x_{2}) = 1 + x_{1} + x_{2}\ .$$
\begin{itemize}
\item[(a)] Check that the manifold $M_P$, corresponding to $P$ under
Delzant's construction, is diffeomorphic to $\cp\#\,\ocp$ (the blow-up of $\cp$
at one point).
\newline \underline{Hint}: read Chapter 1 of~\cite{gui1}.
\item[(b)] Show that the ``canonical" toric K\"ahler metric defined by the
``potential"
$$g_P (x) = \frac{1}{2} \left( \sum_{r=1}^3 \ell_r(x)\log\ell_r(x) +
\ell_{-3}(x)\log\ell_{-3}(x) \right)$$
is not extremal.
\item[(c)] Show that the toric K\"ahler metric defined by the ``potential"
$$g(x_1,x_2) = g_P (x_1,x_2) + \frac{1}{2}\,h(x_1 + x_2)\,,$$
with $h$ a function of one variable satisfying
$$h''(t) = \frac{2}{t^2 + 11 t + 21} - \frac{1}{t+2}\ ,$$
is extremal.
\end{itemize}
\end{exerc}

The toric manifolds $\cp\#\,2\ocp$ and $\cp\#\,3\ocp$ already give 
rise to interesting open questions. $\cp\#\,3\ocp$ is probably the 
simpler and more appealing of the two. When its K\"ahler class is  
$2\pi$ times its first Chern class, the corresponding Delzant polytope
$P\subset\R^2$ is the hexagon defined by
$$\ell_{\pm 1}(x_{1},x_{2}) = 1 \pm x_{1}\,,\ \ell_{\pm 2}(x_{1},x_{2}) = 
1 \pm x_{2}\ \mbox{and}\ \ell_{\pm 3}(x_{1},x_{2}) = 1 \mp 
(x_{1} + x_{2})\ .$$
A simple but somewhat tedious computation shows that the ``canonical''
toric K\"ahler metric defined by the ``potential''
$$g_P (x) = \frac{1}{2} \left( \sum_{r=-3,\,r\neq 
0}^3 \ell_r(x)\log\ell_r(x) \right)$$
is not extremal. It is known, due to an analytic existence result 
first proved by Siu~\cite{siu}, that a K\"ahler-Einstein (hence 
extremal) metric exists and one should be able to 
\begin{center}
{\it explicitly write it down on $P$.} 
\end{center}
Regarding $\cp\#\,2\ocp$, it is known that it does not 
have any K\"ahler metric of constant scalar curvature, but nothing 
prevents it from having extremal K\"ahler metrics. To the author's 
knowledge these have not been constructed or even proved to exist, 
and 
\begin{center}
{\it this question should also be analysed in the set-up presented 
here.}
\end{center}

This symplectic approach to the problem of finding extremal K\"ahler 
metrics on toric manifolds can be considered as an explicit particular 
example of a much more general set-up proposed by Donaldson 
in~\cite{don1} and~\cite{don2}. The part of that set-up more relevant 
for our discussion here is the following. Consider any compact 
symplectic manifold $(M,\om)$, and suppose for simplicity that 
$H^1(M,\R) = 0$. Let $\Jj$ be the space of all complex structures on 
$M$ that are compatible with $\om$, and assume that this space is nonempty.
$\Jj$ can be naturally endowed with the structure of an 
infinite-dimensional K\"ahler manifold. Let $\Gg$ denote the identity 
component of the group of symplectomorphisms of $(M,\om)$, with Lie 
algebra naturally identified with the space $C_{0}^\infty$ of smooth 
functions on $M$ with integral $0$. $\Gg$ acts naturally on $\Jj$, 
and Donaldson shows in~\cite{don1} that an equivariant moment map for 
this action is given by
$$\begin{array}{ccl}
\Jj & \To & \left(C_{0}^\infty\right)^\ast \\
J & \longmapsto & S_{J}\equiv\,\mbox{scalar curvature of}\ 
\om(\cdot,J\cdot)\,,
\end{array}$$
under the pairing
$$(S,H) = \int_{M} S H \frac{\om^n}{n!}\,,\ \forall H\in C_{0}^\infty\,.$$
Hence, critical points of the norm square of the moment map are 
critical points of the square of the $L^2$-norm of the scalar 
curvature, considered now as a functional on the space $\Jj$ of all 
compatible complex structures. In general this space contains several 
different diffeomorphism classes of complex structures, and so the 
critical condition in this context could be more restrictive than the 
extremal condition of Calabi.

When the compact symplectic manifold $(M^{2n}, \om)$ is toric, with 
action $\tau:\T^n\to\Gg$ and moment map $\phi:M\to P\subset\R^n$, we 
can restrict the above considerations to the space $\Jj^{inv}$ 
of invariant compatible complex structures and to the group 
$\Gg^{inv}$ of equivariant symplectomorphisms, with Lie algebra 
naturally identified with the space 
$\left(C_{0}^\infty\right)^{inv}$ of invariant functions on 
$M$ (i.e. smooth functions on $P$) with integral $0$. The moment map 
of the natural action of $\Gg^{inv}$ on $\Jj^{inv}$ is 
again given by the scalar curvature. Since by 
Proposition~\ref{prop:uniqueJ} all complex structures in 
$\Jj^{inv}$ are in the same diffeomorphism class, we do have 
that critical points of the norm square of the moment map are in one 
to one correspondence with extremal K\"ahler metrics. Moreover, 
(\ref{extremalsymp}) gives a symplectic formulation for the complex 
extremal condition saying that the gradient of the scalar curvature is 
an holomorphic vector field:
\begin{itemize}
\item[] an invariant compatible complex structure $J\in 
\Jj^{inv}$ is extremal if and only if the scalar curvature 
$S_{J}$ of the metric $\om(\cdot,J\cdot)$ is a constant plus a linear 
combination of the components $\phi_{1},\,\ldots ,\,\phi_{n}$ of the 
moment map $\phi$.
\end{itemize}

\subsection{Spectral properties of toric manifolds} 
\label{ssec:spectral}

One can define on any Riemannian manifold $M$ the Laplace operator
$$\De : C^\infty (M) \to C^\infty (M)\ .$$
If $M$ is closed, then $\De$ is a self-adjoint elliptic operator on
$L_{1}^2 (M) \equiv$ completion of $C^\infty (M)$ with respect to the 
norm
$$\| \psi \|^2 = \int_{M}\psi^2 + \int_{M} \|\nabla\psi\|^2\ .$$
It follows from the spectral theory of self-adjoint operators that 
$\De$ has discrete eigenvalues
$$0 = \la_{0} < \la_{1}\leq \cdots\leq \la_{j} \leq \cdots,\ 
\la_{j} \stackrel{j\to\infty}{\longrightarrow}\infty \,,$$
and the corresponding eigenfunctions 
$\left\{\psi_{j}\right\}_{j=0}^\infty$, 
satisfying
$$\De\psi_{j} = \la_{j}\psi_{j}, \ \psi_{j}\in C^\infty (M)\cap L_{1}^2 
(M)\,,$$
can be chosen so that they form an orthonormal basis of $L_{1}^2 (M)$.

When one tries to give estimates for the eigenvalues, the Min-Max 
principle plays a fundamental role. It can be formulated as follows. 
Define $H_{j}\subset C^\infty (M),\ j=1,\ldots $, by
$$H_{j} = \left\{ \psi\in C^\infty (M) : \psi\neq 0 \ \mbox{and}\ 
\int_{M}\psi\cdot\psi_{k} = 0,\ k=0,\ldots,j-1 \right\}\ .$$
Then we have that
\begin{equation}\label{min-max}
\la_{j} = {\disps \inf_{\psi\in H_{j}}} \frac{\disps \int_{M} 
\|\nabla\psi\|^2}{\disps \int_{M}\psi^2}\,,\ j= 1, \ldots\ .
\end{equation}

The relation between the spectrum of $\De$ and geometric properties 
of $M$ has a long history, an important chapter of which is Mark Kac's 
celebrated paper~\cite{kac} on the question:
\begin{center}
can one hear the shape of a drum?
\end{center}
When a group acts by isometries on $M$, the Laplace operator can be 
restricted to the subspace of invariant functions and one can 
consider its spectral theory there. The corresponding eigenvalues will 
be called here {\bf invariant eigenvalues}, and one might ask how much 
of the geometry of $M$ can be recovered from them? Our purpose in 
this subsection is to show how asking this question for K\"ahler toric 
manifolds quickly leads to interesting problems.

The Laplace operator $\De$ of a toric K\"ahler metric given in 
holomorphic coordinates $(u,v)$ by~(\ref{metricomp}), has the form
$$
(\De\psi)(u) = -\frac{1}{\det F} \sum_{j,k=1}^n \frac{\p}{\p u_{j}}
\left( (\det F) f^{jk} \frac{\p\psi}{\p u_{k}} \right)
$$
for any $\T^n$-invariant smooth function $\psi = \psi(u)$. 
Using~(\ref{hol-symp}) this can be written in symplectic $(x,y)$ 
coordinates as
\begin{equation}\label{lapsymp}
(\De\psi)(x) = - (\det G)  \sum_{j,k=1}^n g^{jk} \frac{\p}{\p x_{j}}
\left( \frac{1}{\det G}\, \frac{\p\psi}{\p x_{k}} \right)
\end{equation}
where $\psi = \psi(x)$ is again any $\T^n$-invariant smooth function.
Given any Delzant polytope $P$ and K\"ahler toric manifold
$(M_{P},\om_{P},J)$ defined by a 
``potential'' $g$ according to Theorem~\ref{th:Jgen}, the spectrum 
on  $C^\infty(P)$ of the corresponding operator $\De$, 
defined by~(\ref{lapsymp}),  is the invariant spectrum of 
$(M_{P},\om_{P},J)$ and will be denoted by
$$0=\la_{0} < \la_{1}(g) \leq \la_{2}(g) \leq \cdots\ .$$
The corresponding $\T^n$-invariant eigenfunctions are denoted by 
$\psi_{g,j}\in C^\infty(P)$.

One easily computes $\|\nabla\psi\|^2$ to be given pointwise by
$$\|\nabla\psi\|^2 = \sum_{j,k=1}^{n} \frac{\p\psi}{\p x_{j}}
g^{jk} \frac{\p\psi}{\p x_{k}} = G^{-1} \left( \frac{\p\psi}{\p x} 
\right)$$
where $\frac{\p\psi}{\p x} = (\frac{\p\psi}{\p x_{1}}, \ldots , 
\frac{\p\psi}{\p x_{n}})$ and $G^{-1}(\cdot)$ denotes the quadratic 
form on $\R^n$ defined by the symmetric matrix $G^{-1}$. The fact that 
in symplectic coordinates the volume form is always standard, makes 
the dependence of the Min-Max principle~(\ref{min-max}) on the metric 
much easier to understand. In particular, for the invariant 
eigenvalues we get
\begin{equation} \label{evsymp}
\la_{j}(g) = \inf_{\psi\in H_{g,j}} \frac{\disps \int_{P} 
G^{-1} \left( \frac{\p\psi}{\p x} \right)\,dx}
{\disps \int_{p}\psi^2(x)\,dx},\ j=1,\ldots\ ,
\end{equation}
where $H_{g,j} = \left\{ \psi\in C^\infty(P) : \psi\neq 0\ 
\mbox{and}\ \int_{P}\psi\cdot \psi_{g,k} = 0,\ k=0,\ldots,j-1 \right\}$. 
We see that the only major dependence on the metric is 
the $G^{-1}$ that appears on the numerator (for $j=1$ this is 
actually the only dependence on the metric because
$H_{g,1} = \left\{ \psi\in C^\infty(P) : \psi\neq 0\ 
\mbox{and}\ \int_{P}\psi = 0 \right\}$ is independent of $g$). 
Hence, to increase the 
invariant eigenvalues $\la_{j}(g)$ one should try to increase the 
eigenvalues of $G^{-1}$, or equivalently decrease the eigenvalues of 
$G=\Hess_{x}(g)$, always subject to the restrictions imposed by 
Theorem~\ref{th:Jgen}.

In~\cite{abfr} this idea is applied with some interesting results to 
the simplest possible example, the $2$-sphere $S^2$. A theorem of 
J.Hersch~\cite{hers} states that for any smooth metric on $S^2$, with 
total area equal to $4\pi$, the first nonzero eigenvalue of the 
Laplace operator is less than or equal to $2$ (this being the value 
for the standard round metric). It is natural to ask if, for 
$S^1$-invariant metrics, there is a similar upper bound for the first 
invariant eigenvalue? It turns out that the answer in general is no
(see~\cite{eng} and~\cite{abfr}). However one has the following
\begin{theorem}[\cite{abfr}] \label{th:abfr} 
Within the class of smooth $S^1$-invariant metrics $g$ on $S^2$ with total
area $4\pi$ and corresponding to a surface of revolution in
$\R^3$, we have that
$$\la_{j}(g) < \frac{\xi_{j}^2}{2},\ \ j=1,\ldots,$$
where $\xi_{j}$ is the $\left((j+1)/2\right)^{\rm th}$ positive zero of the
Bessel function $J_{0}$ if $j$ is odd, and the $\left(j/2\right)^{\rm th}$ 
positive zero of $J_{0}'$ if $j$ is even. These bounds are optimal.

In particular,
$$\la_{1}(g) < \frac{\xi_{1}^2}{2} \approx 2.89.$$
\end{theorem}
The proof of this theorem goes roughly as follows. Because in dimension $2$ 
any metric is a K\"ahler metric, we have that any $S^1$-invariant metric on 
$S^2$ is a toric K\"ahler metric. 
Hence we can use~(\ref{evsymp}) and try to maximize the single entry 
of the matrix $G^{-1}$ ($n=1$ in this $S^2$ case). It turns out that 
for surfaces of revolution in $\R^3$
one accomplishes that in an optimal way by deforming the standard 
sphere, through a family of ellipsoids of revolution obtained by 
``pressing'' the North and South poles against each other, towards 
the union of two flat discs of area $2\pi$ each (a singular surface).
The estimates of Theorem~\ref{th:abfr} follow from the values of 
the invariant eigenvalues of the Euclidean Laplacian on a disc.

Hersch's theorem has been generalized by J.-P. Bourguignon, P. Li and 
S.T. Yau~\cite{bly} to any K\"ahler metric on a projective complex 
manifold. Hence, it makes sense to generalize the question answered
for $S^1$-invariant metrics on $S^2$ and ask:
\begin{center}
{\it under what conditions does there exist an upper bound for the first 
invariant eigenvalue of a K\"ahler toric manifold?}
\end{center}

Borrowing from M. Kac, we can also ask the following attractive question on
spectral geometry of toric manifolds:
\begin{center}
{\it can one hear the shape of a Delzant polytope?}
\end{center}
As we have seen, associated to every Delzant polytope $P$ we have a 
``canonical'' K\"ahler toric manifold $(M_{P},\om_{P},J_{P})$ given 
by Theorems~\ref{th:delzant} and~\ref{th:Jcan}, and hence a ``canonical'' 
Laplace operator $\De_{P}$, which by~(\ref{lapsymp}) can be explicitly 
written down for invariant functions as
$$\De_{P}(\psi) = - (\det G_{P})  \sum_{j,k=1}^n (g_{P})^{jk}\, 
\frac{\p}{\p x_{j}}\left( \frac{1}{\det G_{P}}\, \frac{\p\psi}{\p x_{k}} 
\right),\ \psi\in C^\infty(P)\,.$$
Its spectrum, which coincides with the invariant spectrum of 
$(M_{P},\om_{P},J_{P})$ and will be denoted by
$$0 = \la_{0} < \la_{1}(P) \leq \la_{2}(P) \leq \cdots\ ,$$
is then a sequence of numbers canonically associated to the Delzant 
polytope $P$, and the question is:
\begin{center}
{\it can we recover $P$ from $\left\{\la_{j}(P)\right\}_{j=1}^\infty$ ?}
\end{center}
Note that, besides being obviously invariant under translations of $P$ 
in $\R^n$, both $\De_{P}$ and its spectrum are also invariant under 
$SL(n,\Z)$ transformations of $\R^n$, i.e.
$$ \left( \De_{P}(f)\right)\circ A^{-1} = \De_{A(P)}(f\circ A^{-1})\,,
\ \forall A\in SL(n,\Z)\,,$$
and so
$$\la_{j}(P) = \la_{j}(A(P))\,,\ j=0,1,\ldots\,,\ 
\forall A\in SL(n,\Z)\,.$$
This follows easily from the fact that a $SL(n,\Z)$ transformation 
of $\R^n$ simply corresponds to a $SL(n,\Z)$ automorphism of the 
acting torus $\T^n$, i.e.
$$(M_{A(P)}, \om_{A(P)}, J_{A(P)}) \cong (M_{P}, \om_{P}, J_{P})\,,$$ 
$$\tau_{A(P)} = \tau_{P}\circ A^T : \T^n \to 
\Diff(M_{A(P)}, \om_{A(P)}, J_{A(P)})\,,$$ 
$$\phi_{A(P)} = A\circ \phi_{P} : M_{A(P)} \to A(P)\subset\R^n\,.$$
Hence, what we can hope to recover from the invariant spectrum
$\left\{\la_{j}(P)\right\}_{j=1}^\infty$ is the equivalence class of
$P$ under translations and $SL(n,\Z)$ transformations of $\R^n$.

\subsection{K\"ahler toric geometry and combinatorics of Delzant polytopes} 
\label{ssec:combinatorics}

The application of geometry of toric varieties to combinatorics of 
convex polytopes has been very successful, with one of the most 
striking examples being the use of the Hard Lefshetz Theorem for 
simplicial toric varieties to prove McMullen's conjectures for the 
number of vertices, edges, faces, etc, of convex simplicial polytopes,
in R.\ Stanley's paper~\cite{stan}. As I.\ Dolgachev writes in his 
review of this paper [MR 81f:52014], ``one must always 
look for toric varieties whenever one has a problem on convex polytopes''.

Our goal in this subsection is to present, in the symplectic setting of
\S\ref{sec:toric}, some tools from K\"ahler geometry that might be 
relevant to combinatorial problems on Delzant polytopes. A particular 
such problem, suggested by R.\ MacPherson during an informal 
conversation in 1997 and directly related to Stanley's result, is the
following.

Let $P\subset\R^n$ be a Delzant polytope and denote by $f_{j}(P),\ 
j=0,\ldots,n,$ the number of $j$-dimensional faces of $P$
(with the convention that $f_{n}(P) = 1$). It is 
well-known (see~\cite{dani} or~\cite{fult}) that the Betti numbers of 
the toric manifold $M_{P}$ associated to $P$ are given by:
\begin{eqnarray}
b^{2k+1}(M_{P}) & = & 0,\ k=0,\ldots,n-1\,; \nonumber \\
b^{2k}(M_{P}) & = & \sum_{j=0}^{k} \left(\matrix n-j \\ n-k 
\endmatrix\right) (-1)^{k-j} f_{n-j}(P),\ k=0,\ldots,n\,. \label{betti}
\end{eqnarray}
The $b^{2k}(M_{P})$ are also denoted by $h^{k}(P),\ k=0,\ldots,n\,,$ 
and called the $h$-numbers of $P$. Because we know that 
$(M_{P},\om_{P})$ is a K\"ahler manifold, the Hard Lefshetz Theorem 
(HLT) applies and says that the map
$$\begin{array}{rrcl}
L^{n-k}\, : & H^{k}(M_{P}) & \to & H^{2n-k}(M_{P}),\ 0\leq k\leq n\,, 
\\ & \al & \mapsto & \al\wedge(\om_{P})^{n-k}
\end{array}$$
is an isomorphism. This immediately puts restrictions on the 
$h$-numbers of $P$, namely:
\begin{itemize}
\item[(i)] $h^{k}(P) = h^{n-k}(P),\ 0\leq k\leq n$; \newline\ 
\item[(ii)] $h^{k+1}(P) - h^{k}(P) \geq 0,\ 
k=0,\ldots,\left[n/2\right]$;
\end{itemize}
which are essentially McMullen's conjectures.

One would like to have a simple combinatorial proof of restrictions 
(i) and (ii) for the $h$-numbers of $P$. Of course, it would be even 
better if one could give a simple combinatorial proof of the full HLT.
On a K\"ahler manifold the simplest proof of HLT is through Hodge 
theory, and the project suggested by MacPherson amounted to
\begin{center}
{\it develop on $P$ a combinatorial version of Hodge theory on 
$(M_{P},\om_{P},J_{P})$.}
\end{center}
The word ``combinatorial'' here means, at least to me, using data and 
analysis coming only from the polytope $P$. Differential forms, 
integration, differentiation, etc, are allowed provided they can be 
defined and performed directly on $P$.

Although MacPherson's project/problem is still very much open, the 
following might be a small initial step towards a solution. As we 
have seen in \S\ref{ssec:can}, associated to a Delzant polytope 
$P\subset\R^n$ we have the affine functions $\ell_{r}:\R^n\to\R$ 
defined by~(\ref{ells}), the potential $g_{P}:P^{\circ}\to\R$ defined 
by~(\ref{gstandard}), its Hessian $G_{P} = \left[ (g_{P})_{jk}\right] =
\Hess(g_{P})$, and the ``canonical'' compatible toric complex structure
$J_{P}$ on $(M_{P},\om_{P})$ given by Theorem~\ref{th:Jcan}. On any 
complex manifold, closed $2$-forms of type $(1,1)$ can be written 
locally as $\p\op$ of a smooth potential function. It is known that the 
cohomology of a K\"ahler toric manifold is generated by its $(1,1)$ 
part, and so one should understand how to represent the $\p\op$ 
operator and the ``generating'' potentials on the polytope $P$. Note 
that, because harmonic forms are always invariant under isometries, 
we only need to understand $\T^n$-invariant forms and potentials
for MacPherson's project.

In holomorphic $(u,v)$ coordinates, and for a $\T^n$-invariant 
function $\nu = \nu(u)$, we have that
$$2i\p\op \nu = \sum_{j,k=1}^n \frac{\disps \p^2 \nu}{\disps \p u_{j} 
\p u_{k}}\, du_{j}\wedge dv_{k}\ .$$
Using~(\ref{hol-symp}) this means that in symplectic $(x,y)$ 
coordinates, and for a $\T^n$-invariant function $\nu = \nu(x)$, we have
\begin{equation}\label{11symp}
2i\p\op \nu = \sum_{j,k=1}^n \frac{\disps \p}{\disps \p x_{j}}
\left( (g_{P})^{kl}\, \frac{\disps \p \nu }{\disps \p x_{l}} \right)\,
dx_{j}\wedge dy_{k}\ .
\end{equation}
\begin{rem} \label{rem:11symp} This formula becomes valid for any 
compatible toric complex structure $J$ on $(M_{P},\om_{P})$, 
determined by a ``potential'' $g\in C^\infty (P^\circ)$ of the form 
given in Theorem~\ref{th:Jgen}, by replacing $g_{P}$ with $g$.
\end{rem}
As to what the relevant ``generating'' potentials should be, Guillemin 
shows in~\cite{gui2} that there is one for each $(n-1)$-dimensional 
face of $P$ and it is given by
\begin{equation}\label{genpot}
\nu_{r}(x) = \log \ell_{r}(x),\ r=1,\ldots,d\equiv f_{n-1}(P)\ .
\end{equation}
He also shows that the differential forms
\begin{equation}\label{genforms}
\al_{r} = \frac{\disps 1}{\disps 2\pi i} \left( \p\op \nu_{r} \right) =
- \frac{\disps 1}{\disps 4\pi} \frac{\disps \p}{\disps \p x_{j}} 
\left( (g_{P})^{kl}\, \frac{\disps \p \log \ell_{r}}{\p x_{l}}\right)\,
dx_{j}\wedge dy_{k}\,,\ r=1,\ldots,d\,,
\end{equation}
represent the Poincar\'e duals to the complex hypersurfaces 
$X_{r}\subset (M_{P},J_{P})$ defined by $\ell_{r}\circ \phi_{P} \equiv 
0$ (i.e. the pre-images under the moment map $\phi_{P}$ of the 
$(n-1)$-dimensional faces of $P$), and are known to generate the full 
cohomology ring of $M_{P}$.
\begin{exerc}[Symplectic potential] \label{ex:stpot}\ \newline
\begin{itemize}
\item[(i)] Using~(\ref{11symp}) and Remark~\ref{rem:11symp} show that,
for any toric complex structure $J$ determined by a ``potential'' 
$g\in C^\infty(P^\circ)$, the function
$$
f_{g}(x) = \sum_{m=1}^{n} x_{m} \frac{\disps \p g}{\p x_{m}}(x) 
\, - \, g(x)
$$
is the potential for the standard symplectic form
$$\om_{P} = \sum_{j=1}^n dx_{j}\wedge dy_{j}\ .$$
Note: this also follows directly from~(\ref{hol-symp}) 
and~(\ref{duality}).
\item[(ii)] Show that 
$$f_{g_{P}}(x) = \frac{\disps 1}{\disps 2} \sum_{r=1}^d \la_{r} \log 
\ell_{r}(x) \, + \, \ell_{\infty}(x)\,,$$
where $\ell_{\infty}(x) = \sum_{r=1}^d \langle x,\mu_{r}\rangle$ is a 
smooth function on the whole $P$.
\item[(iii)] Conclude, using~(\ref{genpot}) and~(\ref{genforms}), that
$$ \frac{\disps 1}{\disps 2\pi} \left[\om_{P}\right] = -
\sum_{r=1}^d \la_{r} \al_{r} \in H^2 (M_{P},\R)\ .$$
\end{itemize}
\end{exerc}

If one applies~(\ref{betti}) to compute $b^2(M_{P})$ one gets
$$b^2(M_{P}) \equiv h^1(P) = f_{n-1} - n \equiv d - n\ .$$
Although~(\ref{genpot}) and~(\ref{genforms}) seem to give more 
generators to $H^2(M_{P})$ than are actually needed, the next 
proposition (stated on purpose in a self-contained and purely combinatorial 
manner) shows that is not the case.
\begin{prop} \label{prop:dim}
Let $P\subset \R^n$ be a Delzant polytope defined by the set of 
inequalities $\ell_{r}(x)\geq 0,\ r=1,\ldots,d$, with
$$\ell_{r}(x) = \langle x,\mu_{r}\rangle - \la_{r}\,,$$
each $\mu_{r}$ being a primitive element of the lattice $\Z^n 
\subset \R^n$ and inward-pointing normal to the $r$-{\rm th} 
$(n-1)$-dimensional face of $P$. Consider the function 
$g_{P}:P^\circ \to \R$ defined by
$$g_{P}(x) = \frac{\disps 1}{\disps 2} \sum_{r=1}^{d} \ell_{r}(x) 
\log \ell_{r}(x)\,,$$
its Hessian $G_{P} = \Hess (g_{P}) = \left[ (g_{P})_{jk} \right]$,
the inverse $G_{P}^{-1} = \left[ (g_{P})^{jk} \right]$, and let 
$\Om^2(P)$ be the real vector space of differential $2$-forms 
generated over $C^\infty(P)$ by
$$dx_{j}\wedge dy_{k},\ 1\leq j,k \leq n\ .$$
Then, the subspace $H^2(P)\subset\Om^2(P)$ generated over $\R$ by the 
differential forms
$$- 4\pi \al_{r} = \frac{\disps \p}{\disps \p x_{j}} 
\left( (g_{P})^{kl}\, \frac{\disps \p \log \ell_{r}}{\p x_{l}}\right)\,
dx_{j}\wedge dy_{k}\,,\ r=1,\ldots,d\,,$$
has real dimension $d-n$.
\end{prop}
\proof{} Consider the linear map $R:\R^d \to H^2(P)$ defined by
$$R(a_{1},\ldots,a_{d}) = \sum_{r=1}^d a_{r} [-4\pi \al_{r}]\ .$$
$R$ is clearly surjective and one easily checks that its kernel 
coincides with the image of the $(d\times n)$ matrix whose rows are 
the coordinates in $\R^n$ of the normals $\mu_{r},\ r=1,\ldots,d$.
Since, by condition (3) of Definition~\ref{def:delzant}, the rank of 
such a matrix is $n$, we conclude that 
$$\dim\,\ker(R) = n\,,$$
and so
$$\dim\,H^2(P) = d-n\,.$$
\QED

Hence we do get a {\bf unique} combinatorial differential form representing 
each cohomology class in $H^2(M_{P},\R)$. Although this has the 
flavour of Hodge theory, the forms $\al_{r}$ are in general not 
harmonic with respect to the ``canonical'' metric 
$\om_{P}(\cdot,J_{P}\,\cdot)$.  For example, the harmonic 
representative for $[\om_{P}] \in H^2(M_{P},\R)$ is always the standard
$$\om_{P} = \sum_{j=1}^n dx_{j}\wedge dy_{j}\,,$$
while we know from Exercise~\ref{ex:stpot} that its representative in 
$H^2(P)$ is given by
$$-2\pi \sum_{r=1}^d \la_{r} \al_{r}\,.$$
It also follows easily from Exercise~\ref{ex:stpot} 
that these two representatives are the same iff
$$\ell_{\infty}(x) = \sum_{r=1}^d \langle x,\mu_{r} \rangle \equiv 
\mbox{constant}\Leftrightarrow \sum_{r=1}^d \mu_{r} = \ov{0} \in\R^n\,,$$
and this is not always true for a Delzant polytope $P$ (take for 
example the $4$-gon or $5$-gon in $\R^2$ corresponding to
$\cp\#\ocp$ and $\cp\# 2\ocp$). Even when this $\ell_{\infty}\equiv 
0$ condition is satisfied, computations on simple examples 
(like $\cp\# 3\ocp$) quickly convince us that the individual forms 
$\al_{r}$ are almost never harmonic (it is likely that the only 
exceptions are projective spaces and their cartesian products).

In spite of this, Proposition~\ref{prop:dim} gives some hope that 
consideration of differential forms of the type given by~(\ref{genforms})
could be an initial step towards a combinatorial Hodge type 
theory on $P$.

\appendix

\section{Proof of Theorem~\ref{th:Jgen}}

In this appendix we give a proof of Theorem~\ref{th:Jgen}. The idea 
is again to translate to symplectic K\"ahler geometry some well-known 
facts from complex K\"ahler geometry.

Let $(M_{P},\om_{P},\tau_{P})$ be the toric symplectic manifold 
associated to a Delzant polytope $P\subset\R^n$, and $J$ any 
compatible toric complex structure. As we know from 
section~\ref{sec:toric}, Delzant's construction also provides 
$(M_{P},\om_{P},\tau_{P})$ with a ``canonical'' compatible toric 
complex structure $J_{P}$. The next proposition explains why, in the 
statement of Theorem~\ref{th:Jgen}, we do not need to impose the 
condition of $J$ being in the same diffeomorphism class as $J_{P}$.
\begin{prop} \label{prop:uniqueJ}
$(M_{P}, J,\tau_{P})$ is equivariantly biholomorphic to
$(M_{P},J_{P},\tau_{P})$.
\end{prop}
\proof{} In the algebraic geometry context this is proved, for 
example, in~\cite{kkms}, Chapter I, \S 2, Theorem 6. The key point of 
that proof is to show that every point of $M_{P}$ admits an open 
invariant affine neighborhood. These neighborhoods are then patched 
together according to the combinatorics of a ``fan'' that uniquely 
determines the resulting normal variety. This same scheme can be 
applied to our symplectic context, with the role of the ``fan'' being 
played by the polytope $P$. We will now sketch how the argument goes.

Every point $p$ of $M_{P}$ admits an open invariant affine 
neighborhood $U_{p}$. The two extreme cases are when $p\in 
M_{P}^\circ$ and $U_{p}$ is equivariantly biholomorphic to $\C^n / 
2\pi i \Z^n$, and when $p$ is a fixed point of the action and 
$(U_{p},p)$ 
is equivariantly biholomorphic to $(\C^n, 0)$ with its standard 
$\T^n$-action. In all cases, the neighborhood itself only depends on 
the action and symplectic structure. The dependence on the complex 
structure appears only in the equivariant biholomorphism. The basic 
reason these neighborhoods exist is that we have $2n$ holomorphic 
vector fields in involution globally defined on $M_{P}$. In fact, 
let $\xi_{1},\ldots,\xi_{n}$ be Hamiltonian vector fields on $M_{P}$ 
induced by the action $\tau_{P}$ from a fixed standard basis of the 
Lie algebra of $\T^n = \R^n / 2\pi\Z^n$. Then 
$(\xi_{1},\ldots,\xi_{n},J\xi_{1},\ldots,J\xi_{n})$ (resp.
$(\xi_{1},\ldots,\xi_{n},J_{P}\xi_{1},\ldots,J_{P}\xi_{n})$) are 
holomorphic vector fields on $(M_{P},J)$ (resp. $(M_{P},J_{P})$) in 
involution, since both $J$ and $J_{P}$ are integrable and $\T^n$-invariant,
and the $\xi_{i}$'s commute with each other. Moreover, on 
$M_{P}^\circ$ these vector fields are all non-zero and linearly 
independent (here the compatibility of $J$ and $J_{P}$ with $\om_{P}$ 
is used) and so can be integrated to give the equivariant 
biholomorphism with $\C^n / 2\pi i \Z^n$. On $M_{P}\setminus 
M_{P}^\circ = \phi_{P}^{-1}(\p P)\,,\ \phi_{P}\equiv$ moment map, the 
vanishing of some of this vector fields is completely determined by 
the combinatorics of $P$, more precisely by the normals to its 
$(n-1)$-dimensional faces, and this information is sufficient to 
construct the required invariant affine neighborhoods of these points.
For example, for a fixed point $p = \phi_{P}^{-1}(\mbox{vertex}\ 
v_{p})$ we have that $U_{p} = \phi_{P}^{-1}(P^\circ \cup \,\{ 
\mbox{faces of $P$ whose closure contains $v_{p}$}\})$ is (both for 
$J$ and $J_{P}$) equivariantly biholomorphic to $(\C^n,0)$.

Having established the existence of the required open invariant affine 
neighborhoods $U_{p}$ for every point $p\in M_{P}$, one then checks 
that the patching of these neighborhoods is also completely 
determined by $P$ and hence is the same for both $(M_{P}, J, 
\tau_{P})$ and $(M_{P}, J_{P}, \tau_{P})$. The simplest example of 
such patching is between $M_{P}^\circ \cong \C^n / 2\pi i \Z^n$ and 
$U_{p}\cong (\C^n , 0)$, for a fixed point $p$ corresponding to a 
vertex $v_{p}$ of $P$ for which the normals to the 
$(n-1)$-dimensional faces that meet at $v_{p}$ are the standard basis 
vectors of $\R^n$. In this case the pathing is simply given by
$$
\begin{array}{rcl}
\C^n / 2\pi i\Z^n & \To & \C^n \\
z = (z_{1},\ldots,z_{n}) & \longmapsto & e^z = (e^{z_{1}}, \ldots,
e^{z_{n}})\,.
\end{array}
$$
\QED

Let $\varphi_{J}:(M_{P},J_{P}) \to (M_{P},J)$ be such an equivariant 
biholomorphism. By the construction given in the above proof,
we know that $\varphi_{J}$ can be chosen so that it acts as the identity 
in cohomology.
Then $(M_{P},\om_{P},J)$ is equivariantly K\"ahler isomorphic to 
$(M_{P},\om_{J},J_{P})$, with $\om_{J} = (\varphi_{J})^\ast (\om_{P})$ 
and $\left[\om_{J}\right] = \left[\om_{P}\right] \in H^2 (M_{P})$.
This means that there exists a $\T^n$-invariant smooth function 
${f_{J}}\in C^\infty (M_{P})$ such that
\begin{equation} \label{a:omhol}
\om_{J} = \om_{P} + 2 i \p\op {f_{J}}\,,
\end{equation}
where the $\p$- and $\op$-operators are defined with respect to the 
complex structure $J_{P}$.

In the $(x,y)$ symplectic coordinates of $M_{P}^\circ \cong 
P^\circ\times\T^n$, obtained via the ``canonical'' moment map 
$\phi_{P}:M_{P}\to\R^n$ with respect to $\om_{P}$, we then have the 
following. By Theorem~\ref{th:Jcan}, $J_{P}$ is obtained from the 
``potential''
$$g_{P}(x) = \frac{1}{2}\sum_{r=1}^{d} \ell_{r}(x) \log \ell_{r}(x)
\in C^\infty (P^\circ)\,,$$
where the $\ell_{r}:\R^n\to\R,\ r=1,\ldots,d$, are the affine 
functions defining $P$ and given by~(\ref{ells}). By~(\ref{11symp}), 
Exercise~\ref{ex:stpot} and~(\ref{a:omhol}), we have that $\om_{J}$ 
is given by
\begin{equation}\label{a:omsymp}
\om_{J} = \om_{P} + 2i\p\op f_{J} = \sum_{j,k=1}^n \frac{\disps \p}{\disps \p x_{j}}
\left( (g_{P})^{kl}\, \frac{\disps \p (f_{P} + f_{J}) }{\disps \p x_{l}} \right)\,
dx_{j}\wedge dy_{k} \equiv 2i\p\op \tilde{f}\,,
\end{equation}
where $\tilde{f}\equiv f_{P} + f_{J}$, $f_{J}\in C^\infty(P)$ and
$$f_{P}(x) \equiv f_{g_{P}}(x) = 
\frac{1}{2}\sum_{r=1}^{d} \ell_{r}(x) \log \ell_{r}(x) + 
\ell_{\infty}(x)\ .$$
Note that, because $\om_{J}(\cdot,J_{P}\,\cdot)$ is a Riemannian metric 
(in particular positive definite), there are restrictions on the 
function $f_{J}$ and these will be specified and used below.

The proof of the first part of Theorem~\ref{th:Jgen} consists now of 
the following three steps:
\begin{itemize}
\item[(i)] write down on $P$ a change of coordinates 
$\tilde{\varphi}_{J}:P\to P$, corresponding to the equivariant 
diffeomorphism $\varphi_{J}:M_{P}\to M_{P}$, that transforms the
symplectic action/angle coordinates $(x,y)$ for $\om_{P}$ into 
symplectic action/angle coordinates $(\tx = 
\tvphi_{J}(x), y)$ for $\om_{J}$;
\item[(ii)] find the ``potential'' $g = g(\tx)$ for the transformed
$J = (\tvphi_{J})_{\ast}(J_{P})$ in this $(\tx,y)$ coordinates;
\item[(iii)] check that the function $h:P^\circ \to\R$, given by
$h(\tx) = g(\tx) - g_{P}(\tx)$ is actually smooth on the whole $P$.
\end{itemize}

The first step is easy. If one examines formula~(\ref{a:omsymp}) 
for $\om_{J}$ attentively, one is immediately led to the change of 
coordinates given in vector form by
$$
\tx = \tvphi_{J}(x) = x + G_{P}^{-1} \cdot \frac{\disps \p f_{J}}{\disps 
\p x}\,,
$$
where $\p f_{J} / \p x = \left( \p f_{J} / \p x_{1}, \ldots , \p f_{J} / 
\p x_{n} \right)^t \equiv$ column vector. Note that the behaviour of the 
matrix $G_{P}^{-1}$ on the boundary of $P$, discussed in 
\S~\ref{ssec:general}, and the fact that $\om_{J}$ given 
by~(\ref{a:omsymp}) is nondegenerate, imply that $\tvphi_{J}$ is a 
diffeomorphism of $P$ that fixes its vertices and preserves each edge,
$2$-face, ..., $n$-face (this last one being the interior $P^\circ$ 
of $P$). The fact that $\om_{J}(\cdot,J_{P}\,\cdot)$ is a Riemannian 
metric implies that
\begin{equation} \label{a:posdef}
\left( d\tvphi_{J}\right) G_{P}^{-1}\ 
\mbox{is symmetric and positive definite on}\ P^\circ\,,
\end{equation}
and so we must also have
\begin{equation} \label{a:posdet}
\det \left( d\tvphi_{J}\right) > 0 \ \mbox{on the whole}\ P
\end{equation}
($d\tvphi_{J}$ denotes the Jacobian matrix of $\tvphi_{J}$). Moreover, 
these conditions also imply that for a point $x_{0}$ belonging to the
$r$-th $(n-1)$-dimensional face of $P$ we have that
\begin{equation} \label{a:positive}
\langle \left( d\tvphi_{J}\right)_{x_{0}}(\mu_{r}),\mu_{r}\rangle > 
0\,,
\end{equation}
where $\mu_{r}$ is the inward pointing normal used in the definition 
of $\ell_{r}$, and $\langle \cdot,\cdot\rangle$ denotes the standard 
euclidean inner product in $\R^n$.

For the second step, the simplest way is to use the Legendre duality 
between the potential for the symplectic form and the ``potential'' 
for the complex structure. In symplectic coordinates this duality is 
given by~(\ref{duality}) (compare also with Exercise~\ref{ex:stpot} 
(i) ). In the $(\tx,y)$ symplectic coordinates for $\om_{J}$ this 
says that
\begin{equation} \label{a:duality}
f(\tx) + g(\tx) = \sum_{m=1}^n \tx_{m} \frac{\disps \p g}{\disps \p 
\tx_{m}}(\tx)\,,
\end{equation}
where $f$ is the $J$-potential for the symplectic form $\om_{J} = 
d\tx\wedge dy$, with the complex structure $J = 
\left(\tvphi_{J}\right)_{\ast}(J_{P})$ being given by 
$$
J = \left[ \begin{array}{ccc}
0 & \vdots & - \left( \Hess_{\tx}(g)\right)^{-1} \\
\cdots & \cdots & \cdots \\
\Hess_{\tx}(g) & \vdots & 0
\end{array}\right]\,.
$$
The following two facts make~(\ref{a:duality}) more explicit for our 
purposes:
\begin{itemize}
\item[(1)] because being the K\"ahler potential of a $(1,1)$-form 
with respect to a complex structure is an invariant notion, we have 
that at $\tx = \tvphi_{J}(x)$
$$f(\tx) = \tilde{f}(x) = f_{P}(x) + f_{J}(x) \,;$$
\item[(2)] because $u = \p g_{P}/\p x\,,\ v=y$ and $\tilde{u} = 
\p g/\p\tx \,,\ \tilde{v} = y$ are holomorphic coordinates for 
$J_{P}$ and $J = \left(\tvphi_{J}\right)_{\ast}(J_{P})$, the map
$\tilde{u} = \tilde{u}(u)\,,\ \tilde{v}=v$ is a biholomorphism. Due 
to its very special form, it can only be the identity plus a constant. 
We will assume, without any loss of generality for our purposes, that 
the constant is zero and so
$$\frac{\disps \p g}{\disps \p\tx_{m}}(\tx) = \left( 
\frac{\disps \p g_{P}}{\disps \p x_{m}}\circ \tvphi_{J}^{-1}\right)
(\tx)\,.$$
\end{itemize}
Using these two facts, we can write~(\ref{a:duality}) as
\begin{equation}\label{a:gtilde}
g (\tx) = \sum_{m=1}^n \tx_{m} \left( \frac{\disps \p g_{P}}
{\disps \p x_{m}}\circ \tvphi_{J}^{-1}\right)(\tx)
- \left( f_{P}\circ \tvphi_{J}^{-1}\right)(\tx)
- \left( f_{J}\circ \tvphi_{J}^{-1}\right)(\tx)\,.
\end{equation}

We can now address the third step. Since $\tvphi_{J}$ is a smooth 
diffeomorphism of $P$, to prove that $h=g -g_{P} \in 
C^\infty (P)$ is equivalent to proving that $h\circ \tvphi_{J}
\in C^\infty (P)$. Using~(\ref{a:gtilde}) and the fact that $f_{J}\in 
C^\infty (P)$, this means that what we need to show is that
the function $\tilde{h}:P^\circ \to\R$ given by
$$
\tilde{h}(x) \equiv h(\tvphi_{J}(x)) + f_{J}(x) = 
\langle\, \tvphi_{J}(x)\,, \frac{\disps \p g_{P}}{\disps 
\p x}(x) \,\rangle - f_{P}(x) - g_{P}(\tvphi_{J}(x))
$$
is smooth on the whole polytope $P$. A simple explicit computation gives
$$ \tilde{h}(x) = \frac{\disps 1}{\disps 2} \left[ 
\ell_{\infty}(\tvphi_{J}(x)) - \ell_{\infty}(x) + \sum_{r=1}^{d}
\ell_{r}(\tvphi_{J}(x)) \log\left(\frac{\disps \ell_{r}(x)}
{\disps \ell_{r}(\tvphi_{J}(x))}\right)\right]\,,$$
and so we will be done if we can prove that the functions
$$
\ga_{r}(x) \equiv \frac{\disps \ell_{r}(x)}
{\disps \ell_{r}(\tvphi_{J}(x))}\,,\ r=1,\ldots,d\,,
$$
are smooth and strictly positive on $P$. The fact 
that they are smooth and strictly positive on $P^\circ$ is immediate. 
The question is, for each $\ga_{r}$, what happens on the $r$-th 
$(n-1)$-dimensional face of $P$? At a point $x_{0}$ belonging to this 
face, and for its normal $\mu_{r}$, we have that
$$d\left( \ell_{r}\circ\tvphi_{J}\right)_{x_{0}}(\mu_{r}) =
d\ell_{r} \left( (d\tvphi_{J})_{x_{0}}(\mu_{r})\right) =
\langle\, \mu_{r} , (d\tvphi_{J})_{x_{0}}(\mu_{r}) \,\rangle > 0$$
by~(\ref{a:positive}). Since $\tvphi_{J}$ preserves the faces of $P$, 
we also have that $\ell_{r}(\tvphi_{J}(x))$ can be written as
$\ell_{r}(\tvphi_{J}(x)) = \ell_{r}(x) \cdot \de_{r}(x)$ for some 
smooth function $\de_{r}$ (at least in a neighborhood of $x_{0}$).
These two together mean that
$$ 0 < d\left( \ell_{r}\circ\tvphi_{J}\right)_{x_{0}}(\mu_{r}) =
\ell_{r}(x_{0}) (d\de_{r})_{x_{0}}(\mu_{r}) + \de_{r}(x_{0}) 
(d\ell_{r})_{x_{0}}(\mu_{r}) = \de_{r}(x_{0}) \cdot \|\mu_{r}\|^2\,,$$
and so, in a neighborhood of $x_{0}$, the function
$$\ga_{r}(x) = \frac{\disps \ell_{r}(x)}
{\disps \ell_{r}(\tvphi_{J}(x))} = \frac{\disps 1}{\disps 
\de_{r}(x)}$$
is smooth and strictly positive as desired.

The facts that $G = \Hess_{\tx}(g)$ is positive definite on 
$P^\circ$ and has determinant of the form
$$\det(G) = \left[ \de (\tx)\cdot\prod_{r=1}^{d} 
\ell_{r}(\tx) \right]^{-1}\ ,$$
with $\de$ a smooth and strictly positive function on the whole 
$P$, follow from~(\ref{a:posdef}) and~(\ref{a:posdet}).

The second part of Theorem~\ref{th:Jgen} can be proved either by 
reversing the reasoning we just did, or directly by showing that any 
complex structure $J$ defined on $M_{P}^\circ \cong P^\circ \times 
\T^n$ by a ``potential'' $g$, satisfying the positivity and 
nondegeneracy conditions specified, compactifies to give a smooth 
compatible toric complex structure on $(M_{P},\om_{P})$.


\subsection*{Acknowledgments}

I want to thank Rosa Dur\~ao for several helpful conversations,
Ana Cannas da Silva for useful comments and suggestions on earlier
versions of this paper, and Yael Karshon for suggesting~\cite{kkms} as a 
possible reference for Proposition~\ref{prop:uniqueJ}.


\end{document}